\documentclass[a4paper,12pt]{article}
    \usepackage[top=2.5cm,bottom=2.5cm,left=2.5cm,right=2.5cm]{geometry}
    \usepackage{cite, amsmath, amssymb}
\usepackage{tikz}
\usepackage{amsmath}

\begin{document}
\begin{center}
{\LARGE\bf No threshold graphs are cospectral}
\end{center}
\begin{center}
{\large \bf J. Lazzarin, O.F. M\'arquez and F. Tura}
\end{center}
\begin{center}
\it Departamento de Matem\'atica, UFSM, Santa Maria, RS, 97105-900, Brazil
\end{center}
\begin{center}
{\small lazzarin@smail.ufsm.br, oscar.f.marquez-sosa@ufsm.br, ftura@smail.ufsm.br}
\end{center}



\newtheorem{Thr}{Theorem}
\newtheorem{Pro}{Proposition}
\newtheorem{Que}{Question}
\newtheorem{Con}{Conjecture}
\newtheorem{Cor}{Corollary}
\newtheorem{Lem}{Lemma}
\newtheorem{Fac}{Fact}
\newtheorem{Ex}{Example}
\newtheorem{Def}{Definition}
\newtheorem{Prop}{Proposition}
\def\floor#1{\left\lfloor{#1}\right\rfloor}

\newenvironment{my_enumerate}{
\begin{enumerate}
  \setlength{\baselineskip}{14pt}
  \setlength{\parskip}{0pt}
  \setlength{\parsep}{0pt}}{\end{enumerate}
}

\newenvironment{my_description}{
\begin{description}
  \setlength{\baselineskip}{14pt}
  \setlength{\parskip}{0pt}
  \setlength{\parsep}{0pt}}{\end{description}
}


\begin{abstract} A threshold graph $G$ on $n$ vertices is defined by binary sequence of length $n.$ In this paper we present an explicit formula for computing the characteristic polynomial of a threshold graph from its binary sequence. Applications include obtaining a formula for the determinant  of adjacency matrix of a threshold graph and showing that no two nonisomorphic  threshold graphs are cospectral. 
\end{abstract}
$\hspace{1cm}${\it keywords:} threshold graph, characteristic polynomial, cospectral graphs.

$\hspace{0.25cm}$ {\it AMS subject classification:} 15A18, 05C50, 05C85.
\baselineskip=0.30in


\section{Introduction}
\label{intro}

Let  $G= (V,E)$ be an undirected graph with vertex set $V$ and edge set $E,$ without loops or multiple edges. 
The {\em adjacency matrix} of $G$, denoted by $A=[a_{ij}]$, is a matrix whose rows and columns are indexed by the vertices of $G$, and is defined to have entries $a_{ij}=1$ if and only if $v_i$ is adjacent to $v_j,$ and $a_{ij}=0$ otherwise.
The characteristic polynomial of $G,$ denoted by $P_{G}(x),$ can be expressed  as
$P_G(x) =\det(A - x I ).$ Their roots  are called eigenvalues of $G,$ or simply the {\em spectrum} of $G.$ 

This paper is concerned  with threshold graphs, introduced by Chv\'atal and Hammer \cite{Chv} and Henderson and Zalcstein \cite{Hen} in 1977. They are an important class of graphs because of their numerous 
applications \cite{Mah95} and have several alternative characterizations. For example, a threshold graph is a $\{2K_2, C_4, P_4\}$-free graph, that is,
a graph which does not contain $2K_2, C_4,$ and $ P_4$ as an induced subgraph. In particular, a threshold graph has a recursive definition based on a binary sequence that will be relevant to this paper, and we will describe in the next Section.

There is a considerable body of knowledge on the spectral properties of threshold graphs
 \cite{chain,Bapat,furer,JTT2015,JTT2014,JTT2013,Stanic,SF2011,Simic}.
However, in the literature, with exception of some particular  cases  does not seem to provide many results about formulas for the characteristic polynomial of graphs. In this paper we attempt to fill this gap, presenting in Section \ref{Secformula}, an explicit  formula to compute the characteristic polynomial of threshold graph from its binary sequence.


 Two nonisomorphic graphs with the same spectrum are called {\em cospectral}. 
In recent years, there has been a growing interest to find families of cospectral graphs. There are many constructions in the literature \cite{godsil82,schwenk71}. This notion is originally defined for the adjacency matrix of the graph $G,$ but a natural extension of the problem is to find families of graphs that are cospectral  with relation to other matrices.
For instance, a family of $2^{n-4}$ threshold graphs on $n\geq 4$ vertices, which are $\mathcal{Q}$-cospectral  was exhibited in  \cite{tura17}, where $\mathcal{Q}$ is the signless Laplacian matrix.

The result in \cite{tura17}  motivates us to investigate further
the cospectrality related  to adjacency matrix. It seemed to
us that cospectrality was present  in the threshold graph class. Indeed, in this
paper we prove that no two threshold graphs are cospectral.
The proof is based on an explicit formula to compute
the characteristic polynomial of threshold graphs. We also use this formula to obtain the determinant of its adjacency matrix.

The paper is organized as follows.
In Section \ref{Sec2}  we show the representation of a threshold graph by binary sequence, and some known results. In Section \ref{Secformula} we present an explicity formula  to compute
the  characteristic polynomial of threshold graphs from its binary sequence. We finalize the paper showing in Section \ref{Sec4} that no two threshold graphs are cospectral.

\section{Preliminaries}
\label{Sec2}
 Recall that  a vertex is  {\em isolated} if it has no neighbors, and is  {\em dominating} if it is  adjacent to all others vertices. 
 A threshold graph is obtained
through an iterative process which starts with an isolated vertex, and where, at each
step, either a new isolated vertex is added, or a dominating vertex is added.
 
 An alternative way to describe the above inductive construction is to associate a threshold graph $G$ of order $n$ 
with a binary  sequence of length $n,$  where a vertex $v_i$ corresponds to digit $0,$  if an isolated  vertex $v_i$ is added and $v_i$  corresponds to digit $1,$ if $v_i$ was added  as a dominating vertex. The choice of digit associated to $v_1$ is arbitrary and we use it as $0.$ Hence, for a connected threshold graph $G$ with $n\geq 2$ vertices the last vertex of $G$ is associated with a digit of type $1.$

\begin{figure}[h!]
       \begin{minipage}[c]{0.25 \linewidth}
\begin{tikzpicture}[ultra thick]
  [scale=0.5,auto=left,every node/.style={circle}]
  \foreach \i/\w in {1/,2/,3/,4/,5/,6/,7/,8/,9/,10/}{
    \node[draw,circle,fill=black,label={360/10 * (\i - 1)+90}:\i] (\i) at ({360/10 * (\i - 1)+90}:3) {\w};} 
  \foreach \from in {3,4,5,9,10}{
    \foreach \to in {1,2,...,\from}
      \draw (\from) -- (\to);}
\end{tikzpicture}
       \end{minipage}\hfill
       \begin{minipage}[l]{0.4 \linewidth}
\begin{tikzpicture}[ultra thick]
\draw[fill=white ](0,0) circle [radius=0.7];
\draw[fill](0, 0.25) circle[radius =0.1];
\draw[fill](0, -0.25) circle[radius =0.1];

\draw[fill=white ](0,3) circle [radius=0.7];
\draw[fill](0, 3.25) circle[radius =0.1];
\draw[fill](-0.25, 2.75) circle[radius =0.1];
\draw[fill](0.25, 2.75) circle[radius =0.1];

\draw(0.7, 3)--(2.3,3);
\draw(3, 2.3)--(3,0.7);
\draw(0.7,0)--(2.3,0);

\draw(0.5, 0.5)--(2.5,2.5);


\draw[fill=gray!30! white](3,3) circle[radius=0.7];
\draw(3, 3.25)--(3,2.75);
\draw[fill](3, 3.25) circle[radius =0.1];
\draw[fill](3, 2.75) circle[radius =0.1];

\draw[fill=gray!30!white ](3,0) circle [radius=0.7];
\draw[fill](3, 0.25) circle[radius =0.1];
\draw[fill](2.75, -0.25) circle[radius =0.1];
\draw[fill](3.25, -0.25) circle[radius =0.1];
\draw(3, 0.25)--(2.75,-0.25);
\draw(3, 0.25)--(3.25,-0.25);
\draw(2.75, -0.25)--(3.25,-0.25);

\end{tikzpicture}
       \end{minipage}
       \caption{A threshold graph and its partitioned representation}
       \label{fig1}
\end{figure}
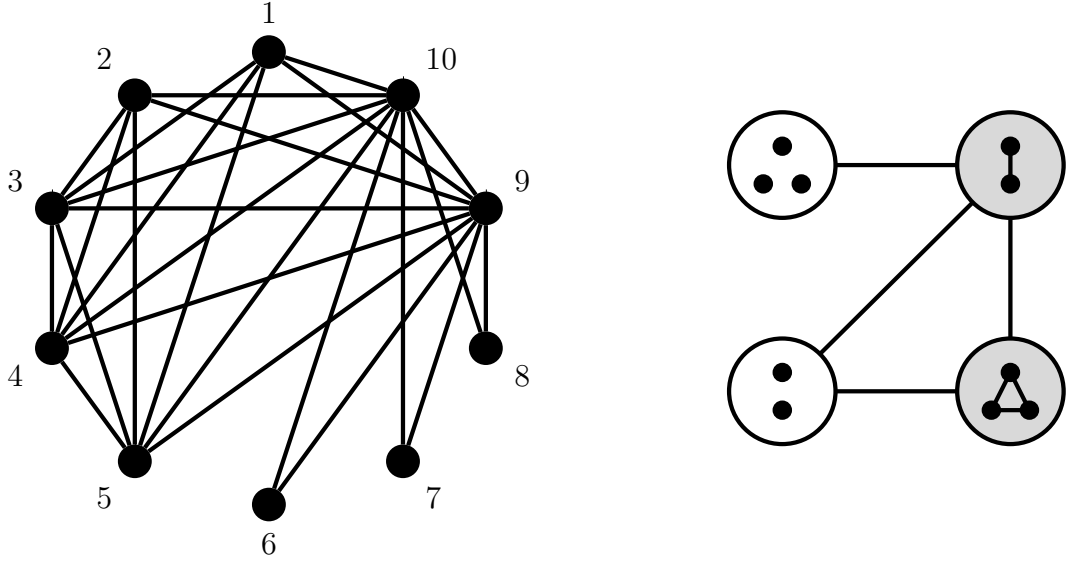

Let  $G= (0^{a_1} 1^{a_2} \ldots 0^{a_{n-1}} 1^{a_n})$  denote a connected threshold graph $G$ 
where each $a_i$ is a positive integer. Let $V_i$ denote the set of vertices corresponding to $a_i $ $ ( i=1,2,\ldots,n).$ Hence, $| V_i | = a_i.$  
It is easy to check that the partition of set vertex of $G$ into non-empty  subsets $V_1, V_2, \ldots, V_n$  determines a divisor $H$ of $G$ \cite{Godsil}. As illustration  the Figure \ref{fig1} shows a threshold graph $G=(0^2 1^3 0^3 1^2)$ and its partitioned representation. The $n\times n$ adjacency matrix $D$ of $H$ 
has the following form

\begin{eqnarray}
\label{eq1}
 D= \left[\begin{array}{ccccccc}
                0   &   a_2            &   0       &  a_4               & \ldots    &0   &  a_{n}  \\
               a_1    &   a_2-1    &   0   &      a_4            & \ldots           & 0& a_{n}\\
               0         &      0        &   0  &     a_4        & \ldots      &     0   &  a_n  \\
               a_1 &  a_2            &    a_3  &  a_4 -1              & \ldots      & 0    & a_{n}    \\
            \vdots         &         \vdots &        \vdots        & \vdots   &  \ddots       &\vdots & \vdots     \\
              0              &           0     &         0      &     \ldots              & 0  & 0& a_n\\
        a_{1}       &    a_{2}     &a_{3}  &\ldots & a_{n-2}  &  a_{n-1}   &     a_n -1\\
\end{array}\right]
\end{eqnarray}

The following theorem (see, for example  \cite{Stanic}) will play an important role in the sequel.

\begin{Thr} 
\label{mainA} For a connected threshold graph $G=(0^{a_1} 1^{a_2} \ldots 0^{a_{n-1}} 1^{a_n}),$  let $\lambda$ be an eigenvalue of $G,$
distinct from $0$ or $-1.$ If $H$ is a divisor of $G$ then $\lambda$ appears    in the spectrum of $H.$
\end{Thr}

We finish this section, presenting an interesting result about the multiplicity of eigenvalues $0$ and $-1$ of a threshold  graph $G.$
Let $m_{0}(G)$ and $m_{-1}(G)$ denote the multiplicity of the eigenvalue $0$ and $-1,$ respectively. 
 According \cite{Bapat,JTT2015},
 they can be obtained  directly from its binary sequence, that is
\begin{Thr}
\label{mainB} For a connected threshold graph  $G= (0^{a_1} 1^{a_2} \ldots 0^{a_{n-1}} 1^{a_n})$
where each $a_i$ is a positive integer. Then   
\begin{itemize}
 \item[i.]  $m_{0}(G) = \sum_{i=1}^{\frac{n}{2}} (a_{2i-1} -1).$ 

\item[ii.] $m_{-1}(G) = \left\{
\begin{array}{lr}
\sum_{i=1}^{\frac{n}{2}} (a_{2i} -1) & \mbox{if $a_1 > 1$   }\\
 1+ \sum_{i=1}^{\frac{n}{2}} (a_{2i} -1) & \mbox{if $a_1=1.$  }
\end{array} \right.$
\end{itemize}
\end{Thr}

\section{An explicity formula for $P_G(x)$ }
\label{Secformula}

 In this section we present an explicit formula for the characteristic polynomial of threshold graph from its binary sequence.
As first application, we obtain an expression to the determinant of adjacency matrix of threshold graph similar to result presented in \cite{Bapat}.

\begin{Thr} 
\label{mainC} Let  $G= (0^{a_1} 1^{a_2} \ldots 0^{a_{n-1}} 1^{a_n})$ be a connected threshold graph
 and let $H$ be divisor of $G.$ For the $n\times n$ adjacency matrix $D$ of $H$ and an indeterminate $x,$ the $D-xI$ matrix can be reduced to tridiagonal matrix of form

\begin{eqnarray}
\label{eq2}
M_{n\times n}=   \left\{\begin{array}{cccc}
           m_{1,1}=x+a_1, m_{i,j} = 0 ,& if  &1< |i-j| &  \\
              m_{k,k} = (-1)^{k+1} a_k ,  &  if    &  1<k \leq n   & \\
          m_{2k+1,2k+2}= -x-1,&  m_{2k,2k +1}= -x&  & \\
          m_{2k+1,2k}= x+1,  &   m_{2k,2k-1}= x & & \\

   \end{array}\right.
 \end{eqnarray}  
\end{Thr}

\noindent {\bf Proof:}
 Let $D$ be the matrix of $H$ given by (\ref{eq1}). We will prove the case $n$ even, the case $n$ odd follows a similar way.  

For the matrix $D-xI,$
  we first perform on the even rows with following operations
$ R_{i}  \leftarrow R_{i} - R_{i-2},$
for each $i =n, n-2, n-4, \ldots,4,$ giving the following  matrix:
\begin{eqnarray}
\label{eq3}
 \left[\begin{array}{cccccccc}
                -x   &   a_2          &   0       &  a_4              & \ldots  &  a_{n-2}& 0 &a_n  \\
               a_1    &   a_2-1-x    &   0   &       a_4           & \ldots         & a_{n-2} & 0  & a_n\\
               0         &      0       &   -x &     a_4        & \ldots   & a_{n-2}  &   0  &     a_n \\
               0 &     x+1        &    a_3 &  -x-1            & \ddots           &     \vdots &   & \vdots    \\
            \vdots         &         \vdots &        \vdots        & \ddots  &  \ddots  &a_{n-2} &  0  & 0    \\
             0       &    0     &   0&x+1 &  a_{n-3}       &- x-1  &  0 &  0 \\
              0              &           0     &         0      &     \ldots        &0  & 0   & -x    & a_n\\
               0       &    0     &   0  &\ldots    &   0  & x+1  &  a_{n-1} &  -x-1 \\
\end{array}\right]
\end{eqnarray}
The next step is to operate on the odd rows of matrix (\ref{eq3}). Then perform the operations
$ R_{i}  \leftarrow R_{i} - R_{i-1},$
for each $i =n -1, n-3, \ldots,3,$ followed by operation $R_2 \leftarrow R_2 -R_1,$ giving the following  matrix:

\begin{eqnarray}
\label{eq4}
 \left[\begin{array}{cccccccc}
                -x   &   a_2          &   0       &  a_4              & \ldots  &  a_{n-2}& 0 &a_n  \\
               x+ a_1    &   -x-1    &   0   &      0           & \ldots         & 0& 0  & 0\\
               x        &      -a_2       &   -x &     0        & \ldots   & 0  &   0  &     0 \\
               0 &     x+1        &    a_3 &  -x-1            & \ddots           &     \vdots &   & \vdots    \\
            \vdots         &         \vdots &        \vdots        & \ddots  &  \ddots  &0 &  0  & 0    \\
             0       &    0     &   0&x+1 &  a_{n-3}       &- x-1  &  0 &  0 \\
              0              &           0     &         0      &     \ldots        &x  & -a_{n-2}   & -x    & 0\\
               0       &    0     &   0  &\ldots    &   0  & x+1  &  a_{n-1} &  -x-1 \\
\end{array}\right]
\end{eqnarray}
For cancel the elements in the first row of matrix (\ref{eq4}), with the exception of the last two elements, we perform the following operations $ R_{1}  \leftarrow R_{1} + R_{i},$
for each $i = 3, 5,\ldots,n-1,$ giving the following  matrix:

\begin{eqnarray}
\label{eq5}
 \left[\begin{array}{cccccccc}
                0  &   0         &   0       &  0             & \ldots  &  0& -x &a_n  \\
               x+ a_1    &   -x-1    &   0   &      0           & \ldots         & 0& 0  & 0\\
               x        &      -a_2       &   -x &     0        & \ldots   & 0  &   0  &     0 \\
               0 &     x+1        &    a_3 &  -x-1            & \ddots           &     \vdots &   & \vdots    \\
            \vdots         &         \vdots &        \vdots        & \ddots  &  \ddots  &0 &  0  & 0    \\
             0       &    0     &   0&x+1 &  a_{n-3}       &- x-1  &  0 &  0 \\
              0              &           0     &         0      &     \ldots        &x  & -a_{n-2}   & -x    & 0\\
               0       &    0     &   0  &\ldots    &   0  & x+1  &  a_{n-1} &  -x-1 \\
\end{array}\right]
\end{eqnarray}
The tridiagonal matrix (\ref{eq2}) is finally obtained after a permutation on the rows of matrix (\ref{eq5}),
followed by a multiplication by $-1$
on the last row. $\hspace{6,5cm} \square$

Let   $Q_n (x) = \prod_{i=1}^n (x -\lambda_i)$ denote the polynomial whose roots are the eigenvalues of a divisor $H$  of $G= (0^{a_1} 1^{a_2} \ldots 0^{a_{n-1}} 1^{a_n}).$
Follows immediately from Theorem \ref{mainC} that the  polynomial $Q_n (x)$  can be computed recursively by
\begin{eqnarray}
\label{eq6}
   \left\{\begin{array}{cc}
           Q_n(x) = (-1)^{n+1}a_n Q_{n-1}(x)+x(x+1)Q_{n-2}(x)& n\geq2  \\
          Q_0(x)=1, &Q_1(x)=x+a_1           \\
   \end{array}\right.
 \end{eqnarray}  
since that the determinant of tridiagonal matrix (\ref{eq2})  may be
easily computed by performing a Laplace expansion.

In order to obtain an explicit formula to $P_G(x)$  for $G= (0^{a_1} 1^{a_2} \ldots 0^{a_{n-1}} 1^{a_n})$ 
we will explore the recursive formula (\ref{eq6}).

Let $[n]=\{1,2,\ldots n\}$, and let $I_{n,l}$ the set of increasing sequences   in $[n]$  of length $l$ such that if $(t_1,t_2,\ldots,t_l)\in I_{n,l}$ then $t_i\equiv n+i-l \, (mod\hspace{0,2cm} 2).$ In other words, elements in $I_{n,l}$ are increasing sequences alternating even and odds numbers such that the last term has the same than parity $n$. For instance
$$I_{7,4}=\{(2,3,4,5),(2,3,4,7),(2,3,6,7), (2,5,6,7),(4,5,6,7) \},$$ while 
$$I_{6,4}=\{(1,2,3,4),(1,2,3,6),(1,2,5,6),(1,4,5,6),(3,4,5,6) \}.$$
In general, for any $l$ the sequences in $I_{7,l}$ must finish in an odd number, while all the sequences in $I_{6,l}$ must finish in an even number.  Given a sequence $\mathbf{t}=(t_1,t_2,\ldots,t_l)$ we denote  $a_\mathbf{t}=a_{t_1} a_{t_2}\cdots a_{t_l}$. Based in this notation we have

\begin{Def}
Let $a_1, a_2, \ldots, a_n$  be fixed sequence of positive integers. We define the following parameter 

  $$   \gamma _{n}(l)=  \left\{\begin{array}{ccc}
           \sum\limits_{\mathbf{t}\in I_{n,l}}  a_\mathbf{t}& if &1\leq l \leq n \\
                1  & if &    l=0      \\
   \end{array}\right.$$
 \end{Def}

\begin{Lem}
\label{lema dos gamas} For $n\geq 4,$
\begin{equation*}
\gamma _{n}(l)=
\begin{cases}
a_{n}\gamma _{n-1}(l-1)+\gamma _{n-2}(l) & \text{ if }1\leq l\leq n-2 \\ 
a_{n}\gamma _{n-1}(l-1) & \text{ if }l\in \{n-1,n\}.
\end{cases}
\end{equation*}
\end{Lem}

\noindent {\bf Proof:}
For $1\leq l\leq n-2$, if $\mathbf{t}=(t_{1},t_{2},\ldots ,t_{l})\in I_{n,l}$
and $t_{l}\neq n$ then the terms in such sequence are less than $n-2$, so then in
this case $\mathbf{t}\in I_{n-2,l}$. Otherwise, $t_{l}=n$ and $1\leq
t_{i}\leq n-1$ satisfies $\mathbf{\hat{t}}=(t_{1},t_{2},\ldots ,t_{l-1})\in
I_{n-1,l-1}$. 
\begin{align*}
\gamma _{n}(l)& =\sum\limits_{\mathbf{t}\in I_{n,l}}a_{\mathbf{t}%
}=a_{n}\sum\limits_{\mathbf{\hat{t}}\in I_{n-1,l-1}}a_{\mathbf{\hat{t}}%
}+\sum\limits_{\mathbf{t}\in I_{n-2,l}}a_{\mathbf{t}} \\
& =a_{n}\gamma _{n-1}(l-1)+\gamma _{n-2}(l).
\end{align*}%
In the case $l\in \{n-1,n\}$ it is clear that $\gamma _{n}(l)=a_{n}\gamma
_{n-1}(l-1)$.  $\hspace{4,5cm} \square$

We write $n=2 m+ r_0$, with $r_0\in \{0,1\}$ and  define $r_1\in\{0,1\}$ such that $r_1\equiv r_0+1 (mod \hspace{0,2cm}2)$. For $n\ge2$ and $y=x+1,$ let 

$$p_n(x,y)=x^{r_0} \sum_{k=0}^{m} (-1)^{m-k}x^k y^k \gamma_n (n-2 k-r_0)  + x^{r_1} \sum_{k=0}^{m-r_1} (-1)^{m-k}x^k y^k \gamma_n(n-2k-r_1).$$
We claim that

\begin{Thr}\label{recorrencia p}
For $n\ge 4$
$$p_n(x,y)=(-1)^{n+1} a_n\,  p_{n-1}(x,y)+  x y\,   p_{n-2}(x,y).$$
\end{Thr}

\noindent {\bf Proof:}
We will prove the case $n=2m$, the case $n=2m+1$ follows in a similar way.
 We can express the polynomial $p$ as:
	$$p_{2m}(x,y)= \sum\limits_{k=0}^{m-1} (-1)^{m-k}x^k y^k\Big( \gamma_{2m} \text{\footnotesize $\big(2(m- k)\big)$}  +x \gamma_{2m}\text{\footnotesize $\big(2(m-k)-1\big)$}\Big)+   x^m y^m .$$
	
	Using the Lemma \ref{lema dos gamas}, we get:
	\begin{align}\label{p2m gamma}
	p_{2m}(x,y)&=(-1)^m (\gamma_{2m}(2m)+x\gamma_{2m}(2m-1))  +\\ \notag
	&\quad  \sum\limits_{k=1}^{m-1} (-1)^{m-k}x^k y^k\Big(a_{2m} \gamma_{2m-1} \text{\footnotesize $\big(2(m- k)-1\big)$}  +a_{2m} x \gamma_{2m-1}\text{\footnotesize $\big(2(m-k)-1-1\big)$}\Big)+\\ \notag
	&\quad  \sum\limits_{k=1}^{m-1} (-1)^{m-k}x^k y^k\Big( \gamma_{2m-2} \text{\footnotesize $\big (2(m-k)\big)$}  +x\gamma_{2m-2}\text{\footnotesize $(2(m-k)-1\big)$}\Big)+   x^m y^m.\\ \notag
	\end{align}
It is a straight computation to see that the sum of the last two terms above
is equal to $xy\,p_{2m-2}(x,y)$. On the other hand, since
\begin{equation*}
(-1)^{m}(\gamma _{2m}(2m)+x\gamma _{2m}(2m-1))=(-1)^{m}(a_{2m}\gamma
_{2m-1}(2m-1)+xa_{2m}\gamma _{2m-1}(2m-2))
\end{equation*}%
We can see that the sum of the first two terms in \eqref{p2m gamma} is $%
(-1)^{2m+1}a_{2m}p_{2m-1}(x,y)$. Thus 
\begin{equation*}
p_{2m}(x,y)=(-1)^{2m+1}a_{2m}p_{2m-1}(x,y)+xy\,p_{2m-2}(x,y)
\end{equation*}
and the proof follows.    $\hspace{11.75cm} \square$

\begin{Cor}
For all $n\ge 2,\quad Q_n(x)=p_n(x,y).$
\end{Cor}

\noindent {\bf Proof:}
	By induction on $n$. If $n=2$, 
\begin{align*}
p_2(x,y)&=x^{0} \sum\limits_{k=0}^{1} (-1)^{1-k}x^k y^k \gamma_2 \text{\footnotesize $(2-2 k-0)$}  +\,  x^{1} \sum\limits_{k=0}^{1-1} (-1)^{1-k}x^k y^k \gamma_2\text{\footnotesize $(2-2k-1)$}\\
&= \left( -\gamma_2 \text{\footnotesize $(2)$} +x y \gamma_2 \text{\footnotesize $(0)$}\right)  -\,  x  \gamma_2\text{\footnotesize $(1)$}\\
&=-a_1a_2+x y-x\, a_2=x y -x a_2 -a_1 a_2=\begin{vmatrix}
x+a_1 & -y \\
x & -a_2 \\
\end{vmatrix}=Q_2(x).
\end{align*}
If $n=3,$
\begin{align*}
p_3(x,y)&=x^{1} \sum\limits_{k=0}^{1} (-1)^{1-k}x^k y^k \gamma_3 \text{\footnotesize $(3-2 k-1)$}  +\,  x^{0} \sum\limits_{k=0}^{1-1} (-1)^{1-k}x^k y^k \gamma_3\text{\footnotesize $(3-2k)$}\\
&= x\left( -\gamma_3 (2) +x y \gamma_3 (0)\right)  +\,   \big(-\gamma_3(3)+x y \gamma_3(1)\big)\\
&=x^2y +xy(a_1+a_3)-x(a_2a_3) -a_1a_2a_3.
\end{align*}
On the other hand
\begin{align*}\begin{vmatrix}
	x+a_1 & -y & 0 \\
	x & -a_2 & -x \\
	0 & y & a_3  \\
\end{vmatrix}&=a_3
\begin{vmatrix}
x+a_1 & -y\\
x & -a_2 
\end{vmatrix}
+x \begin{vmatrix}
x+a_1 & -y  \\
0 & y  \\
\end{vmatrix}\\
&=-a_3a_2(x+a_1)+xy a_3+x(x+a_1)y\\
&=x^2y+xy(a_1+a_3)-a_2a_3x-a_1a_2a_3=Q_3(x).
\end{align*}
By induction hypothesis if $n\ge 4$ and $p_k(x,y)=Q_k(x)$ for all $k<n$, equation (\ref{eq6}) and   Theorem \ref{recorrencia p} conclude $p_n(x,y)=Q_n(x)$ immediately. $\hspace{7cm} \square$

\begin{Thr}
\label{formula}
Let $G= (0^{a_1} 1^{a_2} \ldots 0^{a_{n-1}} 1^{a_n})$ be a connected threshold graph where each $a_i$ is a positive integer.
Let $m_0(G)$ and $m_{-1}(G)$ be the multiplicities of eigenvalues $0$ and $-1$ of $G,$ respectively. Then  $P_G(x)$ is given by 
$$P_G(x) =(-1)^{\sum a_i} x^{m_0(G)} (x+1)^{m_{-1}(G)} Q_n(x), \hspace{0,5cm} where $$ 
$$Q_n(x)=x^{r_0} \sum_{k=0}^{m} (-1)^{m-k}x^k y^k \gamma_n (n-2 k-r_0)  + x^{r_1} \sum_{k=0}^{m-r_1} (-1)^{m-k}x^k y^k \gamma_n(n-2k-r_1),$$
with $r_0, r_1\in\{0,1\},$ such that, $n=2m+r_0, r_1\equiv r_0+1 (mod \hspace{0,2cm}2)$ and $y=x+1.$
\end{Thr}

\begin{Ex}
We apply the formula given in Theorem \ref{formula} to the threshold graph $G=(0^{a_1} 1^{a_2} 0^{a_3} 1^{a_4} 0^{a_5} 1^{a_6})$ with $a_1>1.$ According Theorem \ref{mainA} the multiplicities of $0$ and $-1$ are
$m_0(G) = \sum_{i=1}^{3} (a_{2i-1} -1)$ and $m_{-1}(G) = \sum_{i=1}^{3} (a_{2i} -1).$ The rest of eigenvalues are given by the polynomial $Q_6(x),$ where
\begin{align*}
Q_6(x) &=x^{0} \sum\limits_{k=0}^{3} (-1)^{1-k}x^k y^k \gamma_6 \text{\footnotesize $(6-2 k-0)$}  +\,  x^{1} \sum\limits_{k=0}^{3-1} (-1)^{1-k}x^k y^k \gamma_6\text{\footnotesize $(6-2k-1)$}\\
&=  x^0 \left( -\gamma_6 (6) +x y \gamma_6 (4) - x^2 y^2 \gamma_6(2) +x^3y^3 \gamma_6(0)  \right)  +\\
&  x^{1} \big(-\gamma_6(5)+x y \gamma_6(3)-x^2y^2 \gamma_6(1) \big)\\
&=(-a_1 a_2 a_3 a_4 a_5 a_6  +xy(a_1 a_2 a_3 a_4 + a_1 a_4 a_5 a_6 +a_1 a_2 a_3 a_6) -x^2 y^2(a_1a_2 + \\
& a_1a_4+a_4a_6+a_3 a_4+a_3a_6+a_5a_6) +x^3y^3)-x(a_2a_3a_4a_5a_6)  \\
&  +x^2y (a_2a_3a_4+a_2a_5a_6+a_4a_5a_6) -x^3y^2(a_2 +a_4 +a_6).
\end{align*}
Therefore the $P_G(x)$ of $G=(0^{a_1} 1^{a_2} 0^{a_3} 1^{a_4} 0^{a_5} 1^{a_6})$ is given by
\begin{align*}
P_G(x) &= (-1)^{\sum a_i} x^{m_0(G)} (x+1)^{m_{-1}(G)}   \{x^3 (x+1)^3 -x^3(x+1)^2(a_2+a_4+a_6)    \\
&-x^2 (x+1)^2(a_1a_2 +a_1a_4+a_4a_6+a_3 a_4+a_3a_6+a_5 a_6)  \\
&  +x^2(x+1) (a_2a_3a_4+a_2a_5a_6+a_4a_5a_6) +\\
& x(x+1) (a_1 a_2 a_3 a_4 + a_1 a_4 a_5 a_6 +a_1 a_2 a_3 a_6)\\
& -x(a_2a_3a_4a_5a_6) - a_1a_2a_3a_4a_5a_6\}.
\end{align*}
\end{Ex}

\begin{Cor}
Let $G= (0^{a_1} 1^{a_2} \ldots 0^{a_{n-1}} 1^{a_n})$ be a connected threshold graph where each $a_i$ is a positive integer.
Let $A$ be the adjacency matrix of $G.$ Then $det A =0$ if $a_{2i-1} >1$ for some $i.$ If $a_{2i-1}=1,$ for all $i$ then $$det A =  (-1)^{\sum a_i} \prod_{i=1}^{n} a_{2i}.$$
\end{Cor}
\noindent{\bf Proof:}
Since that the $det(A)$ of $G$  can be computed by $P_G(x)$ valued in $x=0,$
follows immediately from Theorem \ref{formula} that 
$$det(A) =   (-1)^{\sum a_i} 0^{m_0(G)} 1^{m_{-1}(G)} Q_n(0).$$
It is clear that $det(A)=0$ if $a_{2i-1} >1$ for some $i,$ and  $det(A)= (-1)^{\sum a_i} \gamma_n(n) = (-1)^{\sum a_i} \prod_{i=1}^{n} a_{2i}$
if  $a_{2i-1}=1,$ for all $i,$ and the proof follows. $\hspace{5cm} \square$

\section{No threshold graphs are cospectral}
\label{Sec4}
 We finalize this paper showing that no two  nonisomorphic threshold graphs are cospectral. 
For the rest of paper, we let $\gamma _{n}^{\mathbf{a}}{\footnotesize (l)}$ denote the parameter  $\gamma_{n} (l)$
with relation to $\mathbf{a} = (a_1, a_2, \ldots, a_n).$

\begin{Lem}
\label{lem3}
If two threshold graphs  $G= (0^{a_1} 1^{a_2} \ldots 0^{a_{n-1}} 1^{a_n})$  and  $G'= (0^{b_1} 1^{b_2} \ldots 0^{b_{m-1}} 1^{b_m})$ are cospectral, then it must be $n=m.$  
Furthermore, we have  $\gamma _{n}^{\mathbf{a}}{\footnotesize (l)}=\gamma _{n}^{\mathbf{b}}{\footnotesize (l)},$ for  $l=0,1,...,n$.
\end{Lem}
\noindent {\bf Proof:}
It is clear that if $G$ and $G'$ are cospectral they share the same eigenvalues $0$ and $-1,$ with respective multiplicities. Since the remaining eigenvalues  appear in divisors of $G$ and $G',$ respectively,
then follows that  $n=m.$

Let $Q_n(x)$ and $Q'_{n}(x)$ be the polynomials of $G$ and $G',$ from Theorem \ref{formula}.
Since that  $\{1,x,xy,x^{2}y,x^{2}y^{2},x^{3}y^{2},...,x^{k}y^{k}\}$ is
linearly independent set of vectorial space 
$\mathbb{R}\lbrack x],$  follows  from equality  $Q_n(x) = Q'_n(x)$  that 
$\gamma _{n}^{\mathbf{a}}{\footnotesize (l)}=\gamma _{n}^{\mathbf{b}}{\footnotesize (l)},$ for  $l=0,1,...,n.$
$\hspace{3cm} \square$

Let $n=2k$, for $1\le m\le n$ write $m=2l+r_0$ where, $r_0\in\{0,1\}$. Then consider
\begin{equation}
   \Gamma _{n}(m,i)=  
           \sum\limits_{j=0}^{min\{m-1,2i-r_0\}}  (-1)^j  \gamma_n (m-1-j) \gamma_{2i -r_0}(j)
\end{equation}           
We claim that

\begin{Lem} 
\label{formula2}
For $1\leq m \leq n$ 
\begin{equation}
\label{eq9}
\gamma_{n}(m)=\sum\limits_{i=r_0}^{\frac{n-m+r_0}{2}} a_{2i+1-r_0} \Gamma_{n}(m,i)
\end{equation}
\end{Lem}

\noindent{\bf Proof:}
We will prove by induction on $m.$
For $m=1$ since that $r_0 =1,$  then 
\begin{align*}
&  \sum_{i=1}^{k} a_{2i} \Gamma_{n}(1,i) =a_2 \Gamma_{n}(1,1) + a_4  \Gamma_{n}(1,2)  + \ldots + a_{2k}\Gamma_{n}(1,k) \\
& = a_2  \gamma_{2k}(0)\gamma_2(0) + a_4 \gamma_{2k}(0)\gamma_4(0) + \ldots + a_{2k} \gamma_{2k}(0)\gamma_{n}(0) \\
& =  (a_2 + a_4  + a_6 + \ldots  + a_{2k})= \gamma_{n}(1). \end{align*}

We assume that equation (\ref{eq9}) holds for $m$ an even integer. For the case $m$ is an odd integer, the proof is similar.
Since $m+1$ is an odd integer, we have $m+1 = 2l +r_0,$ where $r_0=1.$
We note that $\gamma_{n}(m+1)$ can be viewed as 
\begin{align}\label{eq10}
 &\gamma_n(m+1)= a_2 ( \gamma_{n}(m) - a_1\Gamma_{n}(m,0))  + a_4 (\gamma_{n}(m) -a_1\Gamma_{n}(m,0) - a_{3}\Gamma_{n}(m,1)) \\ \notag
&+ a_6 ( \gamma_{n}(m)  - a_1 \Gamma_{n}(m,0) -a_3\Gamma_n(m,1) - a_{5} \Gamma_{n}(m,2) ) + \ldots   \\ \notag
&+  a_{2k} (\gamma_{n}(m)  - a_1 \Gamma_{n}(m,0) -a_3\Gamma_n(m,1) - a_{5} \Gamma_{n}(m,2) - \ldots -a_{2k-1}\Gamma_{n}(m,k) ) \end{align}


Let $\Gamma'_{n}(m+1,i)$ denote the quantity $\gamma_{n}(m) -\sum_{j=0}^i a_{2j+1}\Gamma_{n}(m,j).$
Then (\ref{eq10}) becomes
\begin{align}\label{eq11}
	\gamma_n(m+1) &= a_2 \Gamma'_{n}(m+1,0)  + a_4 \Gamma'_{n}(m+1,1) + a_6 \Gamma'_{n}(m+1,2) +\\ \notag
	& + a_8 \Gamma'_{n}(m+1,3) + \ldots +  a_{2k}\Gamma'_{n}(m+1,k)\\ \notag
	\end{align}
From equation (\ref{eq11}) it is sufficient to verify  that $$\Gamma'_{n}(m+1,i) = \sum_{j=0}^{min(m,2i-1)}(-1)^j \gamma_{n}(m-j) \gamma_{2i -1}(j).$$
Replacing $\Gamma_{n}(m,i)$ in (\ref{eq10}) and putting the common terms in evidence, we have that 
\begin{align}\label{eq12}
 &\gamma_n(m+1)= a_2 ( \gamma_{n}(m) - a_1\gamma_{n}(m-1))  \\ \notag
&+a_4 (\gamma_{n}(m) -\gamma_{n}(m-1)\gamma_3(1)  +\gamma_{n}(m-2)\gamma_3(2) -\gamma_{n}(m-3)\gamma_3(3)    ) \\ \notag
&+ a_6 ( \gamma_{n}(m)  -\gamma_{n}(m-1)\gamma_5(1) +\gamma_n(m-2)\gamma_{5}(2) + \ldots   \\ \notag
&-\gamma_{n}(m-3)\gamma_5(3)  +\gamma_{n}(m-4)\gamma_5(4) -\gamma_n(m-5)\gamma_{5}(5))+ \ldots   \\ \notag
&+  a_{2k} (\gamma_{n}(m)  -\gamma_{n}(m-1)\gamma_{2k-1}(1) +\gamma_n(m-2)\gamma_{2k-1}(2) \\ \notag
&-\gamma_{n}(m-3)\gamma_{2k-1}(3)+ \ldots   +(-1)^j \gamma_{n}(m-j)\gamma_{2k-1}(j) ) \end{align} 
From (\ref{eq11}) and (\ref{eq12}) the result follows. $\hspace{8.5cm} \square$

\begin{Thr}
\label{cospectral} No two nonisomorphic threshold graphs are cospectral.
\end{Thr}

\noindent{\bf Proof:}
In view of Lemma \ref{lem3} it is sufficient to consider threshold graphs which their respective divisors have the same order. We assume that $G= (0^{a_1} 1^{a_2} \ldots 0^{a_{n-1}} 1^{a_n})$ and
$G'= (0^{b_1} 1^{b_2} \ldots 0^{b_{n-1}} 1^{b_n})$ are cospectral threshold graphs and demonstrate that they must be isomorphic. In other words, we will show that $a_i = b_i,$ for all $i=1,2,\ldots, n.$

By Lemma \ref{lem3} we have that $\gamma_{n}^{\mathbf{a}}(n) =\gamma_{n}^{\mathbf{b}}(n)$ and $\gamma_{n}^{\mathbf{a}}(n-1) =\gamma_{n}^{\mathbf{b}}(n-1),$ that is
\begin{equation}
\label{eq13}
a_1 a_2 \ldots a_{n-1}a_n = b_1 b_2 \ldots b_{n-1}b_n
\end{equation}
and
\begin{equation}
\label{eq14}
 a_2 \ldots a_{n}= b_2 \ldots b_{n}
\end{equation}
By  (\ref{eq13})   and (\ref{eq14})  follows that $a_1= b_1.$ 

By induction, suppose that $a_i=b_i$ for  $i=1,2,\ldots, m-1$ and $m\le n$. We want to prove that $a_m=b_m$. If $m=2 l+r_0$, then $n-m+1=2(k-l)+(1-r_0)$ Using (\ref{eq9}) from Lemma \ref{formula2}, we have:

\begin{align}\label{Gamma1a}
\gamma^{\mathbf{a}}_{n}(n-m+1)&=\sum\limits_{i=1-r_0}^{\frac{m-r_0}{2}} a_{2i+r_0}\,  \Gamma_{n}^{\mathbf{a}}(m,i)\\
&=a_{2-r_0}\Gamma_{n}^{\mathbf{a}}(m,1-r_0)+a_{4-r_0}\Gamma_{n}^{\mathbf{a}}(m,2-r_0)+a_{6-r_0}\Gamma_{n}^{\mathbf{a}}(m,3-r_0)\\ \notag
&+\ldots + a_{m-2} \Gamma_{n}^{\mathbf{a}}(m,\frac{m-r_0}{2}-1)+a_m \Gamma_{n}^{\mathbf{a}}(m,\frac{m-r_0}{2})\notag
\end{align}
Similarly:

\begin{align}\label{Gamma2b}
\gamma^{\mathbf{b}}_{n}(n-m+1)&=\sum\limits_{i=1-r_0}^{\frac{m-r_0}{2}} b_{2i+r_0}\,  \Gamma_n^{\mathbf{b}}(m,i)\\
&=b_{2-r_0}\Gamma_n^{\mathbf{b}}(m,1-r_0)+b_{4-r_0}\Gamma_n^{\mathbf{b}}(m,2-r_0)+b_{6-r_0}\Gamma_n^{\mathbf{b}}(m,3-r_0)\\ \notag
& +\ldots + b_{m-2} \Gamma_n^{\mathbf{b}}(m,\frac{m-r_0}{2}-1)+b_m \Gamma_n^{\mathbf{b}}(m,\frac{m-r_0}{2})\notag
\end{align}
Let $1-r_0\le i\le \frac{m-r_0}{2}$, and let $0\le j\le \min(n-m,2i+r_0-1)$. Since  $\gamma^{\mathbf{a}}_{2i-r_0}(j)$ is a sum and products of $a_1,\ldots a_{m-1}$  we conclude by induction hypothesis that  $\gamma^{\mathbf{a}}_{2i-r_0}(j)=\gamma^{\mathbf{b}}_{2i-r_0}(j)$, for all $j$. Also, we know that  $\gamma^{\mathbf{a}}_{n}(l)=\gamma^{\mathbf{b}}_{n}(l)$, for all $l,$  hence,  $\Gamma_n^{\mathbf{a}}(m,j)=\Gamma_n^{\mathbf{b}}(m,j).$ Thus all the summands in \eqref{Gamma1a} and \eqref{Gamma2b} are equals, except tha last one. Then  

$$a_m \Gamma_n^{\mathbf{a}}(m,\frac{m-r_0}{2})=b_m \Gamma_n^{\mathbf{b}}(m,\frac{m-r_0}{2}). $$

But again, $ \Gamma_n^{\mathbf{a}}(m,\frac{m-r_0}{2})= \Gamma_n^{\mathbf{b}}(m,\frac{m-r_0}{2}) $, and we conclude that $a_m=b_m.$ $\hspace{1.5cm} \square$

\end{document}